\documentclass{article}
\usepackage{lipsum}
\usepackage{layout}
\usepackage{amsmath}
\usepackage{graphicx}
\usepackage{times}
\usepackage{xcolor}
\usepackage{amsmath}
\usepackage{mathtools}
\usepackage{amssymb}
\usepackage{epstopdf}
\usepackage[numbers,sort&compress]{natbib}
\usepackage{appendix}
\newcommand{\be}{\begin{equation}
\newcommand{\ee}{\end{equation}}}
\newcommand{\bea}{\begin{eqnarray}}
\newcommand{\eea}{\end{eqnarray}}
\newcommand{\nn}{\nonumber}
\begin{document}

\title{\textbf{The solution of conformable Laguerre  differential equation using conformable Laplace transform}}
\author{Eqab.M.Rabei, Ahmed Al-Jamel, Mohamed.Al-Masaeed \\
Physics Department, Faculty of Science,\\ Al al-Bayt University, P.O. Box 130040, Mafraq 25113, Jordan\\
eqabrabei@gmail.com, eqabrabei@aabu.edu.jo\\
aaljamel@aabu.edu.jo, aaljamel@gmail.com \\ moh.almssaeed@gmail.com
}

\maketitle


\begin{abstract}
In this paper,  the conformable  Laguerre  and associated Laguerre differential equations are solved using the Laplace transform. The solution is found to be in exact agreement with that obtained using the power series. In addition some of properties of the Laguerre polynomial is discussed and  the conformable Rodriguez's Formula and generating function are proposed and proved
\\
\end{abstract}
\maketitle
\section{Introduction}
Special functions have developed from a wide array of practical challenges that attract not only mathematicians, but also other academics in science who want to learn more about their qualities, characteristics, and applications. The special functions are induced as solutions of well-known differential equations, one of which is the Laguerre differential equation. \\
A specific sort of integral transform is the Laplace Transform. When a function $f (t)$ is considered, the appropriate Laplace Transform is $\mathcal{L}[f(t)]$, where $\mathcal{L}$ is the operator applied to the time domain function $f(t)$. A function's Laplace Transform yields a new function with complex frequency s. The Laplace Transform, like the Fourier Transform, is used to solve differential and integral equations. It's also widely used in the investigation of transient occurrences in electrical circuits, where frequency domain analysis is employed \cite{asher2013introduction}.
The fractional derivative has played an important role in physics,mathematics and engineering  sciences \cite{oldham1974fractional,miller1993introduction,kilbas2006theory,klimek2002lagrangean,agrawal2002formulation,baleanu2006fractional,rabei2004potentials,rabei2006quantization,rabei2007hamilton,rabeihamilton}.The definition of fractional derivative and fractional integral subject to several approaches\cite{oldham1974fractional,miller1993introduction,kilbas2006theory,klimek2002lagrangean,agrawal2002formulation,baleanu2006fractional,rabei2004potentials,rabei2006quantization,rabei2007hamilton,rabeihamilton,hilfer2000applications,podlubny1998fractional} i.e:Riemann-Lioville fractional derivatives, Caputo,Riesz,
Riesz-Caputo,Weyl,Grünwald-Letnikov,Hadamard,
and chen derivatives.\\
Recently, a new definition of derivative of fractional  order was presented by Khalil et.al \cite{khalil2014new} called the conformable derivative. This definition is a natural extension of the usual derivative.The conformable derivative of the constant is zero, and satisfies the standard properties of the traditional derivative i.e the derivative of the product and the derivative of the quotient of two functions and satisfies the chain rule.In addition, one can say the conformable derivative is simple and singular to the standard derivative.\\
For this reason,the last few years, the conformable fractional calculus is applied successfully in various fields (math \cite{atraoui2021existence,atangana2015new,singh2018numerical,shihab2021associated,khalil2019geometric,hammad2014abel,ahmad2015antiperiodic}, physics \cite{al2021extension,chung2021new,AlMasaeedRabeiAlJamelBaleanu+2021+395+401,al2021wkb,al2021effect,mozaffari2018investigation,alextension,hammad2021analytical}, modeling \cite{yavuz2018conformable,xin2019modeling,kumar2020variety}, biology \cite{kumar2018modified,khater2019modified})  .Besides, The conformable fractional Euler-Lagrange equation and Hamiltonian formulation were discussed by Lazo and Torres\cite{lazo2016variational} and the laplace transform in conformable derivative is presented by Abdeljawad \cite{abdeljawad2015conformable}.

The power series and the Laplace transformation are two different methods for solving differential equations \cite{boyce2017elementary}. The conformable laguerre differential equation is solved using conformable power series in ref \cite{abu2019laguerre,shat2019fractional}. Recently, even though Abdelhakim \cite{abdelhakim2019flaw}, pointed that the conformable derivative is not fractional derivative but it is a derivative of fractional order and it is simple and applicable which makes it of interests . Thus, we call it conformable derivative. In  this paper, we want to solve the conformable Laguerre and  associated Laguerre differential equation using conformable Laplace transform  and discuss  their properties. Besides, the Conformable Rodriguez's Formula is  proposed and proved. 
\section{Conformable Laplace transform}
we start by presenting some definitions related to our work.\\
\textbf{Definition 2.1.} we denote a function $f\in [0,\infty) \to R$.The conformable  derivative of $f$ with order $\alpha$ is defined by \cite{khalil2014new}
\be
\label{conformable}
T_\alpha(f)(t)=\lim_{\epsilon \to 0}\frac{f(t+\epsilon t^{1-\alpha})-f(t)}{\epsilon}
\ee
\\
\textbf{Definition 2.2.}.$I_\alpha^a (f)(t)=I_1^a (t^{\alpha-1}f)= \int_a^t \frac{f(x)}{x^{1-\alpha}} dx$ where the integral is the usual Riemann improper integral and $\alpha \in (0,1)$.
\\
\textbf{Definition 2.3.} let $t_0 \in R , 0<\alpha \leq 1$ and $f:[t_0 ,\infty) \to $ be real valued function.Then the conformable  Laplace transform of order $\alpha$ is defined by 
\be
\label{def}
\mathcal{L}_\alpha^{t_0}=F_\alpha^{t_0} (s)=\int_{t_0}^\infty e^{-s\frac{(t-t_0)^\alpha}{\alpha}} f(t)d\alpha(t,t_0)=\int_{t_0}^\infty e^{-s\frac{(t-t_0)^\alpha}{\alpha}} f(t)(t-t_0)^{\alpha-1}dt.
\ee
\\
Following to ref \cite{abdeljawad2015conformable}, the conformable Laplace transform for some certain functions 
\be \label{ex1}(i)\quad \mathcal{L}_\alpha^{t_0} [1]=\frac{1}{s} , s>0\ee

\be \label{ex2}(ii)\quad\mathcal{L}_\alpha^{t_0} [t]=\frac{t_0}{s} +\alpha^{\frac{1}{\alpha}} \frac{\Gamma(1+\frac{1}{\alpha})}{s^{1+\frac{1}{\alpha}}} , s>0\ee

\be \label{ex3}(iii)\quad \mathcal{L}_\alpha^0 [t^p]=\alpha^{\frac{p}{\alpha}} \frac{\Gamma(1+\frac{p}{\alpha})}{s^{1+\frac{p}{\alpha}}} , s>0\ee

\be \label{ex4}(iv)\quad\mathcal{L}_\alpha^0 [e^{\frac{t^\alpha}{\alpha}}]=\frac{1}{s-1} , s>1\ee

\be \label{ex5}(v)\quad\mathcal{L}_\alpha^0 [\sin{\omega \frac{t^\alpha}{\alpha}}]=\frac{1}{w^2+s^2} \ee

\be \label{ex6}(vi)\quad\mathcal{L}_\alpha^0 [\cos{\omega \frac{t^\alpha}{\alpha}}]=\frac{s}{w^2+s^2}\ee
\\
\textbf{Some properties}:let $f$ and $g : \in [0,\infty$ and let $\lambda , \mu , a \in R $ and $0<\alpha \leq 1$. then \cite{al2019new,khader2017conformable}

\be 
\label{pro1} 
(i)\quad
\mathcal{L}_\alpha [\lambda  f+\mu g](t)=(\lambda \mathcal{L}_\alpha f+\mu \mathcal{L}_\alpha g)(t)=\lambda F_\alpha(s)+ \mu G_\alpha (s) , s>0
\ee

\be 
\label{pro2} 
(ii)\quad
\mathcal{L}_\alpha[e^{-a\frac{t^\alpha}{\alpha}} f(t)]=F_\alpha(s+a) , s>|a|
\ee
\be 
\label{pro3} 
(iii)\quad
\mathcal{L}_\alpha[\frac{t^{n\alpha}}{\alpha^n} f(t)]=(-1)^n \frac{d^n}{ds^n} F_\alpha (s), s>0
\ee
\be 
\label{pro4} 
(iv)\quad
\mathcal{L}_\alpha[f*g](t)=F_\alpha(s)*G_\alpha(s)
\ee

\be 
\label{pro5} 
(v)\quad
\mathcal{L}_\alpha[T_\alpha f(t)]=sF_\alpha(s)-f(0)
\ee
\be 
\label{pro6} 
(v)\quad
\mathcal{L}_\alpha[T_{2\alpha} f(t)]=s^2F_\alpha(s)-sf(0) \quad, \quad 0<\alpha\leq \frac{1}{2}
\ee
\\
 In this paper, we adopt $D^\alpha f$ to denote the conformable  derivative of $f$ with order $\alpha$.
\section{Conformable Laguerre  differential equation}
The conformable laguerre  differential equation can be written by replacing the integer derivative by the conformable derivative of order $\alpha$ as 
\be
\label{laguerre fractional}
x^\alpha D^\alpha D^\alpha y +(\alpha-x^\alpha)D^\alpha y + n \alpha y=0
\ee
where $0<\alpha \leq 1$ and $n$ is natural number.\\
The conformable Laplace transform of this equation is
\be
\label{eq1}
\mathcal{L}_\alpha[x^\alpha D^\alpha D^\alpha y]+\alpha \mathcal{L}_\alpha[D^\alpha y]-\mathcal{L}_\alpha[x^\alpha D^\alpha y]+ n \alpha \mathcal{L}_\alpha[y]=0.
\ee
Making use of conformable Laplace properties we have
\bea
\nn
  \mathcal{L}_\alpha[x^\alpha D^\alpha D^\alpha y]&=&-\alpha \frac{d}{ds}\mathcal{L}_\alpha[ D^\alpha D^\alpha y]=-\alpha \frac{d}{ds}[s^2 Y_\alpha(s)-s y(0)]\\\nn 
  &=&-\alpha s^2 \frac{d}{ds}Y_\alpha(s)-2s\alpha  Y_\alpha(s)+\alpha y(0),
  \\ \label{first term} 
  &=& -\alpha s^2 \frac{d}{ds}Y_\alpha(s)-2s\alpha  Y_\alpha(s)+\alpha
\eea
where $y(0)=1 $ and

\bea
\nn
  \mathcal{L}_\alpha[D^\alpha y]&=&s Y_\alpha(s) -y(0)  \\ \label{second term} &=&s Y_\alpha(s) -1.
\eea
Also,
\bea
\nn
     \mathcal{L}_\alpha[x^\alpha D^\alpha y]&=&-\alpha \frac{d}{ds}\mathcal{L}_\alpha[  D^\alpha y]=-\alpha \frac{d}{ds}[s Y_\alpha(s) -y(0)]\\ \nn 
     &=&-\alpha[s\frac{d}{ds} Y_\alpha(s)+Y_\alpha(s)].
 \eea

Then, Eq.\eqref{eq1} takes the form
\be
\label{diff}
-\alpha s(s-1)\frac{d}{ds} Y_\alpha(s)+\alpha( n +1-s)Y_\alpha(s)=0,
\ee
where $\mathcal{L}_\alpha[y]=Y_\alpha(s)$ and divideding  by $-\alpha s(s-1)$, we have
\be
\label{diff2}
\frac{d}{ds} Y_\alpha(s)-\frac{(n +1-s)}{ s(s-1)}Y_\alpha(s)=0.
\ee
Then, the solution read as 
\be
\label{diff slo}
Y_\alpha(s)=\frac{(s-1)^n}{s^{n +1}}.
\ee
 Which can be written as
 \be
\label{diff slo1}
Y_\alpha(s)=\frac{1}{s} (1-\frac{1}{s})^n.
\ee
 It can be expanded in series representation as following, 
 \be
 \label{diff repre}
  Y_\alpha(s)=\sum_{k=0}^n (-1)^k \frac{n!}{(n-k)!  k!}\frac{1}{s^{k+1}}.
 \ee
 Using inverse conformable Laplace transform $\mathcal{L}_\alpha^{-1}$, we have
 \be
 \label{diff inverse}
  \mathcal{L}_\alpha^{-1} Y_\alpha(s)=y(x)=\sum_{k=0}^n (-1)^k \frac{n!}{(n-k)!  k!} \mathcal{L}_\alpha^{-1}\frac{1}{s^{k+1}}. 
 \ee
Setting $k=\frac{p}{\alpha}$, then
\bea
\nn
   \mathcal{L}_\alpha^{-1}\frac{1}{s^{\frac{p}{\alpha}+1}}&=&\frac{\alpha^{\frac{p}{\alpha}}\Gamma(1+\frac{p}{\alpha})}{\alpha^{\frac{p}{\alpha}}\Gamma(1+\frac{p}{\alpha})}\mathcal{L}_\alpha^{-1}\frac{1}{s^{\frac{p}{\alpha}+1}}, \\\nn 
   &=&\frac{1}{\alpha^{\frac{p}{\alpha}}\Gamma(1+\frac{p}{\alpha})}\mathcal{L}_\alpha^{-1}\frac{\alpha^{\frac{p}{\alpha}}\Gamma(1+\frac{p}{\alpha})}{s^{\frac{p}{\alpha}+1}},\\\nn &=&\frac{1}{\alpha^{\frac{p}{\alpha}}\Gamma(1+\frac{p}{\alpha})}x^p.
\eea
Now, replacing  $p$ by $ k \alpha$, we get 
\be
\label{inverse}
\mathcal{L}_\alpha^{-1}\frac{1}{s^{k+1}}=\frac{1}{\alpha^k \Gamma(1+k)}x^{k \alpha}=\frac{1}{\alpha^k k!}x^{k \alpha}.
\ee
Thus, the solution of Eq.\eqref{laguerre fractional} takes the form
\be
\label{sol laguerre fractional}
y(x)= \sum_{k=0}^n (-1)^k \frac{n!}{\alpha^k (n-k)!  (k!)^2} x^{k \alpha}.
\ee
 Here, we obtained the solution of  the laguerre conformable  differential equation using conformable Laplace transform and it is found to be in exact agreement with the solution obtained  using conformable power series, in Ref \cite{abu2019laguerre,shat2019fractional}   .
\subsection{ Conformable Rodriguez's Formula}
In this subsection we propose the conformable Rodriguez formula in the following form
\be
\label{main rod formula}
  L_n(\frac{ x^\alpha}{\alpha})= \frac{e^{\frac{x^\alpha}{\alpha}}}{\alpha^n  n!} D^{n\alpha}[x^{n\alpha}e^{-\frac{x^\alpha}{\alpha}}].
\ee
Which is equivalent to conformable Laguerre  polynomial.\\ 
\textbf{Proof}. Using conformable Leibniz rule 
\be
\label{Leibinz rule}
D^{n\alpha}(f g)= \sum_{k=0}^n\left( \begin{array}{c}
      n\\
k      
\end{array}\right) D^{(n-k)\alpha} f D^{k\alpha} g.
\ee
we get
\bea
\nn
   D^{n\alpha}[x^{n\alpha}e^{-\frac{x^\alpha}{\alpha}}]= \sum_{k=0}^n\left( \begin{array}{c}
      n\\
k      
\end{array}\right) D^{(n-k)\alpha} x^{n\alpha} D^{k\alpha} e^{-\frac{x^\alpha}{\alpha}},
\eea
\be
\label{l1}
 D^{n\alpha}[x^{n\alpha}e^{-\frac{x^\alpha}{\alpha}}]=\sum_{k=0}^n \frac{(-1)^k \alpha^{n-k}(n!)^2}{(k!)^2(n-k)!} x^{k\alpha} e^{-\frac{x^\alpha}{\alpha}}.
\ee
Now,substituting in eq.\eqref{main rod formula}
\bea
\nn
    L_n(\frac{ x^\alpha}{\alpha})&=& \frac{e^{\frac{x^\alpha}{\alpha}}}{\alpha^n  n!}\sum_{k=0}^n \frac{(-1)^k \alpha^{n-k}(n!)^2}{(k!)^2(n-k)!} x^{k\alpha} e^{-\frac{x^\alpha}{\alpha}},\\\nn
    &=&\sum_{k=0}^n \frac{(-1)^k n!}{\alpha^k (k!)^2(n-k)!} x^{k\alpha}.
\eea
So, our proposed formula is equivalent to the conformable polynomial eq.\eqref{sol laguerre fractional}.\\
Figures 1 to 5 depict the behavior of conformable Laguerre functions of several orders for different values of $\alpha$. One can see that as $\alpha$ approaches 1, the Laguerre functions approaches to the traditional one  

\begin{figure}[h!]
    \centering
    \includegraphics[width=0.68\textwidth]{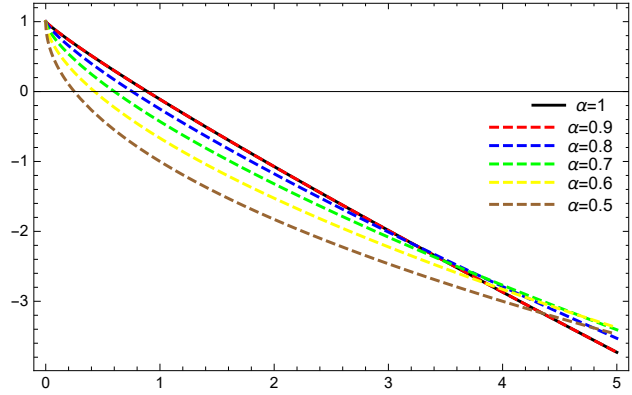}
    \caption{Plot of $ L_1(\frac{ x^\alpha}{\alpha})=1-\frac{x^\alpha}{\alpha}$}
    \label{fig:my_label}
\end{figure}
\begin{figure}[h!]
    \centering
    \includegraphics[width=0.68\textwidth]{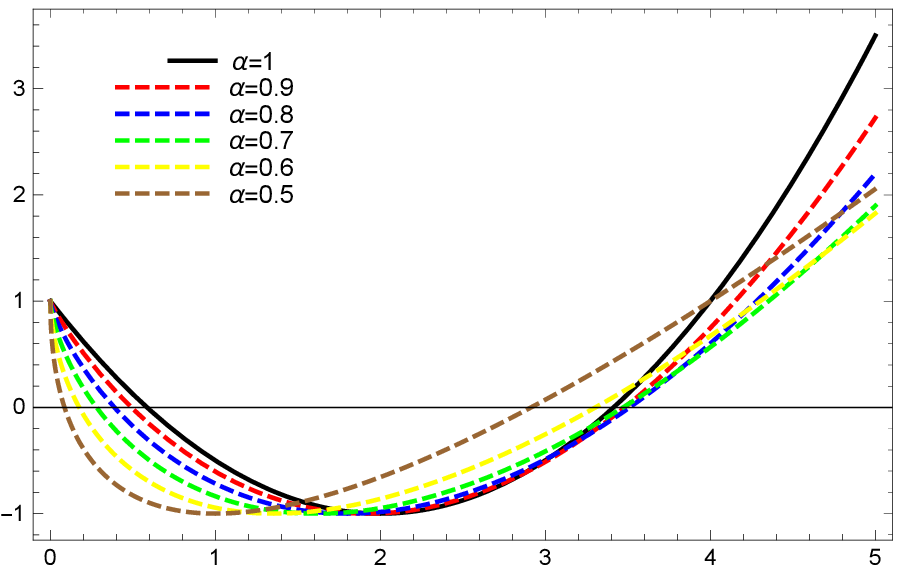}
    \caption{Plot of $ L_2(\frac{ x^\alpha}{\alpha})=1+\frac{x^{2\alpha}}{2\alpha^2}-2\frac{x^\alpha}{\alpha}$}
    \label{fig:my_label}
\end{figure}
\newpage
\begin{figure}[h!]
    \centering
    \includegraphics[width=0.68\textwidth]{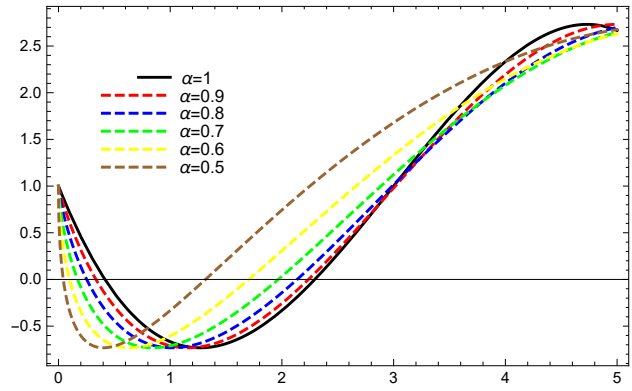}
    \caption{ Plot of$L_3(\frac{ x^\alpha}{\alpha})=-\frac{x^{3\alpha}-9 \alpha x^{2\alpha}+ 18 \alpha^2 x^{\alpha}- 6\alpha^3}{6\alpha^3} $}
    \label{fig:my_label}
\end{figure}
\begin{figure}[h!]
    \centering
    \includegraphics[width=0.68\textwidth]{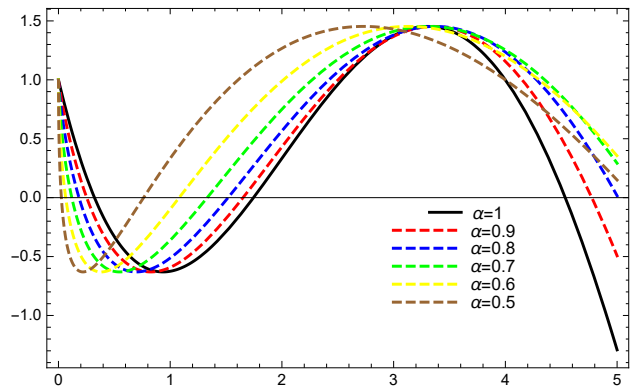}
    \caption{  Plot of $L_4(\frac{ x^\alpha}{\alpha})= 1+\frac{x^{4\alpha}}{24\alpha^4}-2\frac{x^{3\alpha}}{3\alpha^3}+ 3 \frac{x^{2\alpha}}{\alpha^2}-4\frac{x^\alpha}{\alpha}$}
    \label{fig:my_label}
\end{figure}
\newpage
\begin{figure}[h!]
    \centering
    \includegraphics[width=0.68\textwidth]{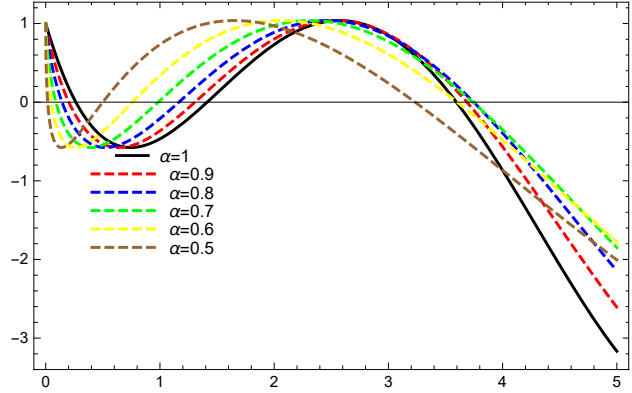}
    \caption{  Plot of $L_5(\frac{ x^\alpha}{\alpha})= 1-\frac{x^{5\alpha}}{120\alpha^5}+5\frac{x^{4\alpha}}{24\alpha^4}-5\frac{x^{3\alpha}}{3\alpha^3}+ 5 \frac{x^{2\alpha}}{\alpha^2}-5\frac{x^\alpha}{\alpha}$}
    \label{fig:my_label}
\end{figure}

The conformable Laguerre polynomials developed here allowed us to gradually manage the polynomials we employed by adjusting the alpha value. We can see how conformable Laguerre polynomials behave when alpha values change without affecting their behavior in these figures.
\subsection{Some remarks  of conformable Laguerre   functions}
\textbf{1.} The values of the  conformable Laguerre  polynomial at $x=0$, read as
\bea
\label{pro1}
&(i)& L_n(0)=1 , \\\label{pro2} &(ii)& D^\alpha L_n(0)= -k, \\\label{pro3}  &(iii)& D^\alpha D^\alpha L_n(0) = \frac{n(n-1)}{2}.
\eea
\textbf{Proof.} $(i)$Following to Ref.\cite{abu2019laguerre}, the generating function of the conformable laguerre  functions read as 
\be
\label{gene fun}
\frac{1}{1-t} e^{-\frac{x^\alpha t}{\alpha(1-t)}}= \sum_{n=0}^\infty L_n(\frac{x^\alpha}{\alpha}) t^n, \quad | t | <1 
\ee
Setting $x=0$, we recoverd
\bea
\nn
\frac{1}{1-t} = \sum_{n=0}^\infty L_n(0) t^n.
\eea
Then, we obtain $L_n(0)=1$ where  $\frac{1}{1-t} = \sum_{n=0}^\infty t^n$. \\$(ii)$ The general Laguerre conformable differential equation read as 
\bea
\nn
x^\alpha D^\alpha D^\alpha L_n(x^\alpha) +(\alpha-x^\alpha)D^\alpha L_n(x^\alpha) + k \alpha L_n(x^\alpha)=0,
\eea
taking $x=0$, we have 
\bea
\nn
\alpha D^\alpha L_n(0) + k \alpha L_n(0)=0.
\eea
But $L_n(0) =1$, then  $D^\alpha L_n(0)=-k$.\\ $(iii)$ Taking the conformable derivative of order $\alpha$ of both sides with respect to 
 eq.(\ref{gene fun}) twice, we have the following relation
\bea
\nn
\frac{1}{1-t} \frac{ t^2}{(1-t)^2} e^{-\frac{x^\alpha t}{\alpha(1-t)}}= \sum_{n=0}^\infty D^\alpha D^\alpha L_n(\frac{x^\alpha}{\alpha}) t^n,
\eea
substituting $x=0$, we get 
\bea
\nn
 \frac{ t^2}{(1-t)^3}= \sum_{n=0}^\infty D^\alpha D^\alpha L_n(0) t^n,
\eea
the left hand side can be expanded as 
\bea
\nn
t^2 [1+3t+\frac{3\cdot 4 }{2!} t^2 + \frac{3\cdot4\cdot5 }{3!} t^3+ \dots + \frac{3\cdot4\cdot5\cdot \dots \cdot n }{(n-2)!}t^{n-2}+\dots] = \sum_{n=0}^\infty D^\alpha D^\alpha L_n(0) t^n.
\eea
Thus, Equating coefficient of $t^n$, we get 
\bea
\nn
D^\alpha D^\alpha L_n(0) = \frac{3\cdot4\cdot5\cdot \dots \cdot n }{(n-2)!}=\frac{ n! }{(n-2)!} =  \frac{n(n-1)}{2}.
\eea
\textbf{2.} The conformable Laguerre  functions are orthogonal 
\bea
\label{ortho rela}
\int_0^\infty e^{-\frac{s^\alpha}{\alpha}}   L_n(\frac{ x^\alpha}{\alpha})  L_m(\frac{ x^\alpha}{\alpha}) d^\alpha x = \delta_{nm} 
\eea
\textbf{ Proof.} Using the generating function 
\bea
\nn
\frac{1}{1-t} e^{-\frac{x^\alpha t}{\alpha(1-t)}}&=& \sum_{n=0}^\infty L_n(\frac{x^\alpha}{\alpha}) t^n, 
\eea
and replacing $t$ by $s$
\bea
\nn 
\frac{1}{1-s} e^{-\frac{x^\alpha s}{\alpha(1-s)}}&=& \sum_{m=0}^\infty L_m(\frac{x^\alpha}{\alpha}) s^m.
\eea
Thus, we have 
\bea
\nn
\frac{1}{(1-t)(1-s)} e^{-\frac{x^\alpha t}{\alpha(1-t)}} e^{-\frac{x^\alpha s}{\alpha(1-s)}} = \sum_{n=0}^\infty  \sum_{m=0}^\infty L_n(\frac{x^\alpha}{\alpha}) L_m(\frac{x^\alpha}{\alpha}) s^m t^n.
\eea
Multiplying both sides by $e^{-\frac{x^\alpha}{\alpha}}$ and and integrating over $d^\alpha x$, we have
\bea
\nn
\frac{1}{(1-t)(1-s)} \int_0^\infty e^{-\frac{x^\alpha}{\alpha}}  e^{-\frac{x^\alpha t}{\alpha(1-t)}} e^{-\frac{x^\alpha s}{\alpha(1-s)}} d^\alpha x = \sum_{n=0}^\infty  \sum_{m=0}^\infty \int_0^\infty e^{-\frac{x^\alpha}{\alpha}} L_n(\frac{x^\alpha}{\alpha}) L_m(\frac{x^\alpha}{\alpha}) s^m t^n d^\alpha x.
\eea
The integral in the left hand side can be calculated as 
\bea
\nn
\int_0^\infty e^{-\frac{x^\alpha}{\alpha}[1+\frac{t}{(1-t)}+\frac{s}{(1-s)}]}  x^{\alpha-1} dx = \frac{1}{[1+\frac{t}{(1-t)}+\frac{s}{(1-s)}]}.
\eea
Substituting this result, we have 
\bea
\nn
\frac{1}{(1-st)} = \sum_{n=0}^\infty  \sum_{m=0}^\infty \int_0^\infty e^{-\frac{x^\alpha}{\alpha}} L_n(\frac{x^\alpha}{\alpha}) L_m(\frac{x^\alpha}{\alpha}) s^m t^n d^\alpha x.
\eea
Then we have the value of the integral as series representations as
\bea
\nn
\frac{1}{(1-st)} = \sum_{n=0}^\infty s^n t^n.
\eea
Then the general relation becomes
\bea
\nn
\sum_{n=0}^\infty s^n t^n = \sum_{n=0}^\infty  \sum_{m=0}^\infty \int_0^\infty e^{-\frac{x^\alpha}{\alpha}} L_n(\frac{x^\alpha}{\alpha}) L_m(\frac{x^\alpha}{\alpha}) s^m t^n d^\alpha x.
\eea
It is clear that, if $m=n $, we have
\bea
\nn
\int_0^\infty e^{-\frac{x^\alpha}{\alpha}} L_n(\frac{x^\alpha}{\alpha}) L_m(\frac{x^\alpha}{\alpha})  d^\alpha x = 1.
\eea 
and if $m \neq n $, we have 
\bea
\nn
\int_0^\infty e^{-\frac{x^\alpha}{\alpha}} L_n(\frac{x^\alpha}{\alpha}) L_m(\frac{x^\alpha}{\alpha})  d^\alpha x = 0.
\eea 
\section{Conformable associated Laguerre  differential equation}
The conformable associated Laguerre  differential equation is defined by 
\be
\label{Associatedlaguerre fractional}
x^\alpha D^\alpha D^\alpha y +(m \alpha+\alpha-x^\alpha)D^\alpha y + n \alpha y=0
\ee
\textbf{Theorem.} If $S(x)$ is a solution of conformable Laguerre  differential equation of order $(m+n)$, then $D^{m\alpha} S $ will be a solution of  conformable associated Laguerre  differential equation.\\
\textbf{Proof.} The conformable Laguerre differential equation of order $(m+n)$ can be written as
\be
\label{laguerre fractional nm}
x^\alpha D^\alpha D^\alpha S +(\alpha-x^\alpha)D^\alpha S + (n+m) \alpha S =0.
\ee
Taking conformable derivative m times $D^{m\alpha}$, we get 
\bea
\nn
D^{m\alpha} [x^\alpha D^\alpha D^\alpha S ] + D^{m\alpha} [(\alpha-x^\alpha)D^\alpha S] + D^{m\alpha} [ (n+m) \alpha S]=0
\eea
Making use of  Leibniz rule, we get 
\bea
\nn
 x^\alpha D^{(m+2)\alpha} S + m \alpha D^{(m+1)\alpha} S +  \alpha D^{(m+1)\alpha} S -x^\alpha D^{(m+1)\alpha} S - m \alpha D^{m\alpha} +  (n+m) \alpha D^{m\alpha}S =0,
\eea
which can be written as 
\bea
\nn
 x^\alpha D^{(m+2)\alpha} S + (m \alpha  +  \alpha  - x^\alpha ) D^{(m+1)\alpha} S  +  n \alpha D^{m\alpha}S =0.
\eea
This equation can be rearranged as
\bea
\nn
 x^\alpha D^{\alpha} D^{\alpha}(D^{m\alpha}S) + (m \alpha  +  \alpha  - x^\alpha ) D^{\alpha} (D^{m\alpha}S)  +  n \alpha (D^{m\alpha}S) =0.
\eea
Thus, if $S(x)$ is a solution of conformable Laguerre differential equation, then $D^{m\alpha} S$ is a solution of the conformable associated Laguerre differential equation.\\
So, we obtain $L_n^m (\frac{x^\alpha}{\alpha}) = D^{m\alpha} L_{n+m} (\frac{x^\alpha}{\alpha})$, then from the definition, we get 
\be
\label{Dm}
L_n^m (\frac{x^\alpha}{\alpha}) =(-1)^m D^{m\alpha} L_{n+m} (\frac{x^\alpha}{\alpha}).
\ee
This is the conformable associated Laguerre  functions. 
\subsection{ Conformable Associated Laguerre  functions}
The conformable associated Laguerre  polynomial can be written as  
\be
\label{frac asso poly}
L_n^m (\frac{x^\alpha}{\alpha})=\sum_{r=0}^{n} (-1)^{r} \frac{(n+m)!}{\alpha^{r} (n-r)!  (r+m)! r!}  x^{r \alpha}.
\ee
\textbf{Proof.} According  eq.(\ref{sol laguerre fractional}), we can write the conformable  Laguerre  functions of order $(m+n)$ as 
\be
L_{m+n} (\frac{x^\alpha}{\alpha}) = \sum_{k=0}^n (-1)^k \frac{(n+m)!}{\alpha^k (n+m-k)!  (k!)^2} x^{k \alpha}.
\ee
Substituting it in eq.(\ref{Dm}), we get 
\bea
\nn
L_n^m (\frac{x^\alpha}{\alpha})&=& (-1)^m D^{m\alpha} \sum_{k=0}^n (-1)^k \frac{(n+m)!}{\alpha^k (n+m-k)!  (k!)^2} x^{k \alpha},\\\nn
&=&   \sum_{k=m}^{n+m} (-1)^{k+m} \frac{(n+m)!}{\alpha^k (n+m-k)!  (k!)^2} \frac{k! \alpha^m }{(k-m)!} x^{(k-m) \alpha}, \\\nn
&=&   \sum_{k=m}^{n+m} (-1)^{k+m} \frac{(n+m)!}{\alpha^{k-m} (n+m-k)!  k! (k-m)!}  x^{(k-m) \alpha}.
\eea
Put $r=k-m$, we have the eq.(\ref{frac asso poly}). 
One may recover the formula (\ref{frac asso poly}) by setting $r=k-m$
\subsection{Conformable Rodriguez's formula}
In this subsection, we show that the conformable Rodriguez formula  for the conformable  associated Laguerre  polynomial takes the form 
\be
\label{rod formula}
  L_n^m(\frac{ x^\alpha}{\alpha})= \frac{x^{- m \alpha }e^{\frac{x^\alpha}{\alpha}}}{\alpha^n  n!} D^{n\alpha}[x^{(n+m)\alpha}e^{-\frac{x^\alpha}{\alpha}}].
\ee
\textbf{Proof}. Using Leibniz rule 
\be
\label{Leibinz rule}
D^{n\alpha}(f g)= \sum_{r=0}^n\left( \begin{array}{c}
      n\\
r      
\end{array}\right) D^{(n-r)\alpha} f D^{r\alpha} g
\ee
we get
\bea
\nn
   D^{n\alpha}[x^{(n+m)\alpha}e^{-\frac{x^\alpha}{\alpha}}]= \sum_{k=0}^n\left( \begin{array}{c}
      n\\
r     
\end{array}\right) D^{(n-r)\alpha} x^{(n+m)\alpha} D^{r\alpha} e^{-\frac{x^\alpha}{\alpha}}.
\eea
Which can be written as, 
\be
\label{l1}
 D^{n\alpha}[x^{(n+m)\alpha}e^{-\frac{x^\alpha}{\alpha}}]=\sum_{r=0}^n \frac{(-1)^r \alpha^{n-r}(n+m)! n!}{(r+m)!(n-r)! r!} x^{(r+m)\alpha} e^{-\frac{x^\alpha}{\alpha}}.
\ee
Now,substituting in eq.\eqref{rod formula}
\bea
\nn
    L_n^m(\frac{ x^\alpha}{\alpha})&=& \frac{x^{- m \alpha }e^{\frac{x^\alpha}{\alpha}}}{\alpha^n  n!}\sum_{r=0}^n \frac{(-1)^r \alpha^{n-r}(n+m)! n!}{(r+m)!(n-r)! r!} x^{(r+m)\alpha} e^{-\frac{x^\alpha}{\alpha}},\\\nn
    &=&\sum_{r=0}^n \frac{(-1)^r (n+m)! }{\alpha^{r} (r+m)!(n-r)! r!} x^{r\alpha} .
\eea
It is equal to eq.(\ref{frac asso poly}).\\ Figures 6 to 11 depict the behavior of conformable  associated Laguerre functions of different orders for some values of $\alpha$.
\begin{figure}[h!]
    \centering
    \includegraphics[width=0.68\textwidth]{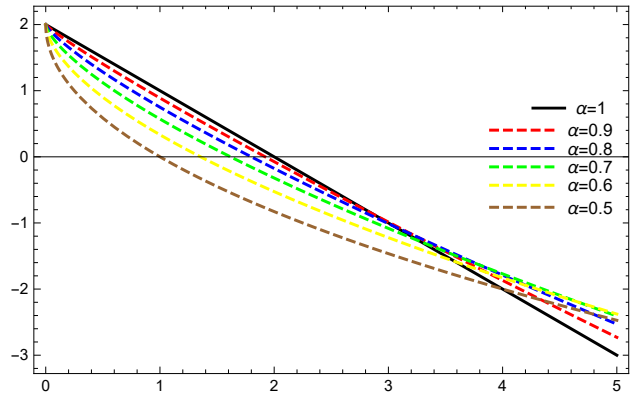}
    \caption{Plot of $ L_1^1(\frac{ x^\alpha}{\alpha})=2-\frac{x^\alpha}{\alpha}$}
    \label{fig:my_label}
\end{figure}
\newpage
\begin{figure}[h!]
    \centering
    \includegraphics[width=0.68\textwidth]{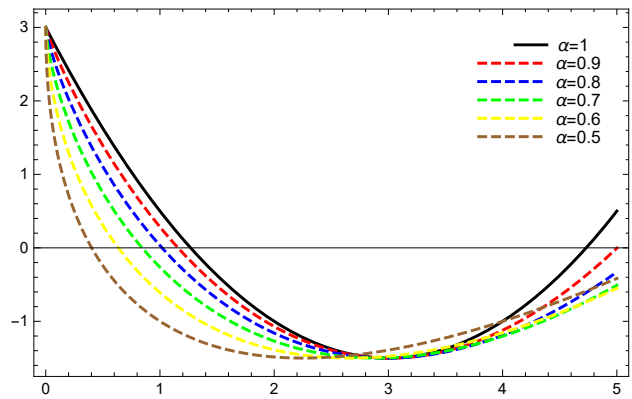}
    \caption{Plot of $ L_2^1(\frac{ x^\alpha}{\alpha})=\frac{x^{2\alpha}}{2\alpha^2}-3\frac{x^\alpha}{\alpha}+3$}
    \label{fig:my_label}
\end{figure}
\begin{figure}[h!]
    \centering
    \includegraphics[width=0.68\textwidth]{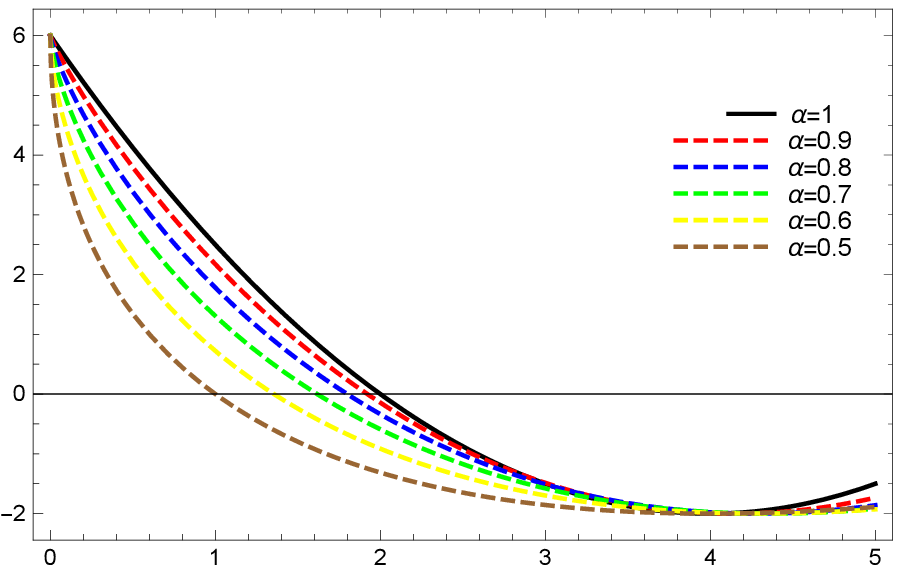}
    \caption{ Plot of$L_2^2(\frac{ x^\alpha}{\alpha})=\frac{x^{2\alpha}}{2\alpha^2}-4\frac{x^\alpha}{\alpha}+6 $}
    \label{fig:my_label}
\end{figure}
\newpage
\begin{figure}[h!]
    \centering
    \includegraphics[width=0.68\textwidth]{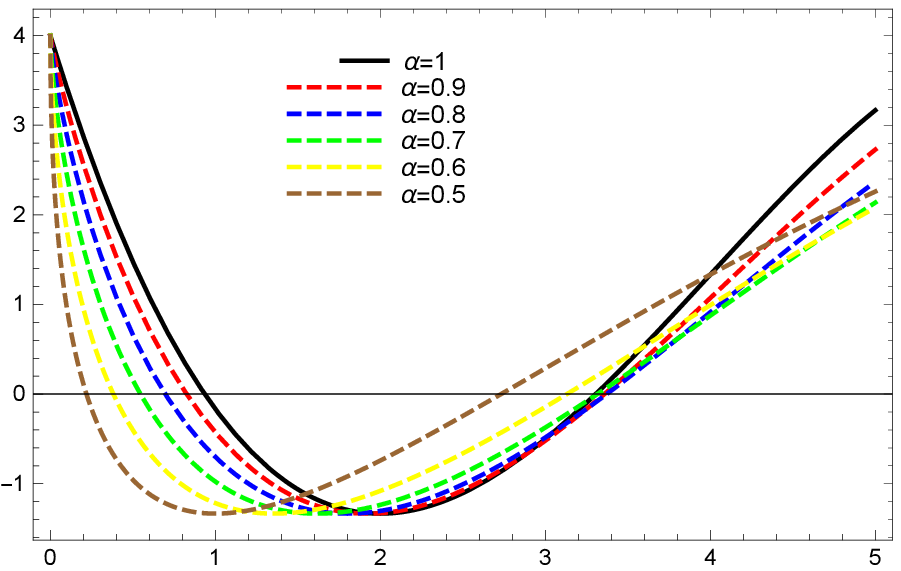}
    \caption{  Plot of $L_3^1(\frac{ x^\alpha}{\alpha})= -\frac{x^{3\alpha}}{6\alpha^3}+ 2 \frac{x^{2\alpha}}{\alpha^2}-6\frac{x^\alpha}{\alpha}+4$}
    \label{fig:my_label}
\end{figure}
\begin{figure}[h!]
    \centering
    \includegraphics[width=0.68\textwidth]{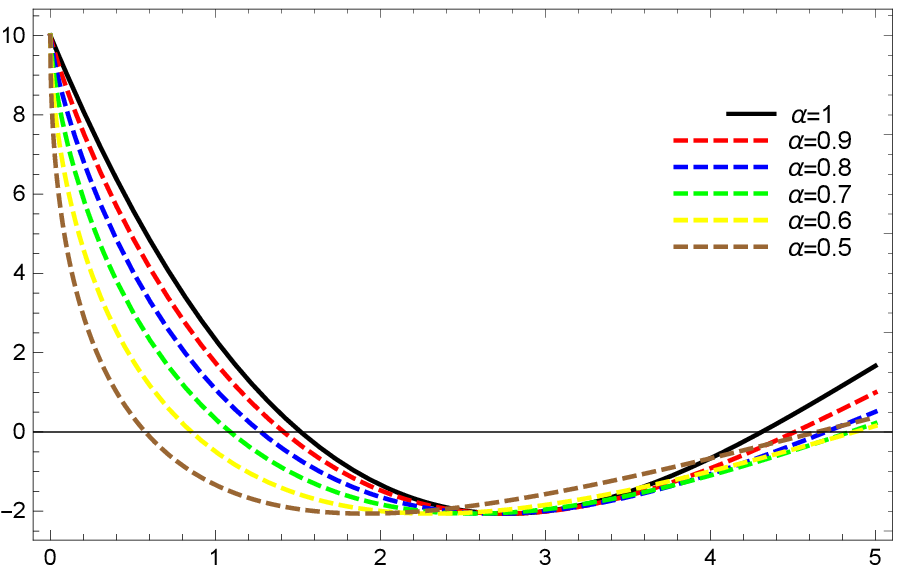}
    \caption{  Plot of $L_3^2(\frac{ x^\alpha}{\alpha})= -\frac{x^{3\alpha}}{6\alpha^3}+ 15 \frac{x^{2\alpha}}{6 \alpha^2}-10\frac{x^\alpha}{\alpha} + 10$}
    \label{fig:my_label}
\end{figure}
\newpage
\begin{figure}[h!]
    \centering
    \includegraphics[width=0.68\textwidth]{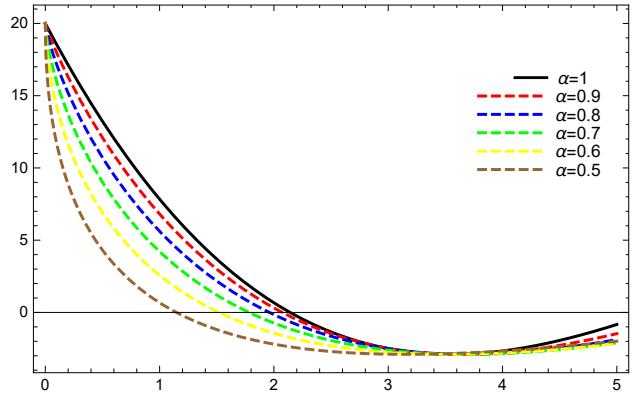}
    \caption{  Plot of $L_3^3(\frac{ x^\alpha}{\alpha})= -\frac{x^{3\alpha}}{6\alpha^3}+ 3 \frac{x^{2\alpha}}{ \alpha^2}-15\frac{x^\alpha}{\alpha} + 20$}
    \label{fig:my_label}
\end{figure}
 The conformable associated Laguerre functions developed here allowed us to gradually manage the polynomials we employed by adjusting the alpha value. We can see how conformable associated Laguerre polynomials behave when alpha values change without affecting their behavior in these figures.

\subsection{Generating function }
The generating function for the conformable associated Laguerre  polynomial is proposed as
\be
\label{gene fun associ}
\frac{ e^{-\frac{x^\alpha t}{\alpha (1-t)}}}{(1-t)^{m+1}} = \sum_{n=0}^\infty L_n^m(\frac{ x^\alpha}{\alpha}) t^n.
\ee
Using eq.(\ref{gene fun}) and  differentiating it $m$ times, we have 
\bea
\nn
 D^{m\alpha} \frac{1}{1-t} e^{-\frac{x^\alpha t}{\alpha(1-t)}}=  D^{m\alpha} \sum_{n=0}^\infty L_n(\frac{x^\alpha}{\alpha}) t^n, \\\nn 
 \frac{(-1)^m t^m}{(1-t)^{m+1}} e^{-\frac{x^\alpha t}{\alpha(1-t)}}=  D^{m\alpha} \sum_{n=m}^\infty L_n(\frac{x^\alpha}{\alpha}) t^n.
\eea
In the R.H.S put $n=m+r$ , we have 
\bea
\nn
 \frac{(-1)^m t^m}{(1-t)^{m+1}} e^{-\frac{x^\alpha t}{\alpha(1-t)}}=  D^{m\alpha} \sum_{r=0}^\infty L_{m+r}(\frac{x^\alpha}{\alpha}) t^{m+r},
\eea
from eq.(\ref{Dm}), we get 
\bea
\nn
 \frac{(-1)^m t^m}{(1-t)^{m+1}} e^{-\frac{x^\alpha t}{\alpha(1-t)}}=   \sum_{n=0}^\infty  (-1)^m  L_n^m (\frac{x^\alpha}{\alpha}) t^{m+n}.
\eea
It is equal the eq.(\ref{gene fun associ}). \\\\
\section{Conclusions}
We have solved the conformable Laguerre and associated Laguerre differential equations using the Laplace transform, we have observed the result found to be in exact agreement with the solution obtained using conformable power series, in Ref \cite{abu2019laguerre,shat2019fractional}. In addition, we have given good proposal for Rodriguez's formula and the generating function and we shown that the conformable Laguerre functions are orthogonal. Moreover from the figures $1-11$ one may observe that, the obtained results coincide with those of regular   Laguerre and associated Laguerre functions as $\alpha$ approaches 1. 
\bibliography{ref}

\begin{thebibliography}{10}
\providecommand{\url}[1]{#1}
\csname url@samestyle\endcsname
\providecommand{\newblock}{\relax}
\providecommand{\bibinfo}[2]{#2}
\providecommand{\BIBentrySTDinterwordspacing}{\spaceskip=0pt\relax}
\providecommand{\BIBentryALTinterwordstretchfactor}{4}
\providecommand{\BIBentryALTinterwordspacing}{\spaceskip=\fontdimen2\font plus
\BIBentryALTinterwordstretchfactor\fontdimen3\font minus
  \fontdimen4\font\relax}
\providecommand{\BIBforeignlanguage}[2]{{%
\expandafter\ifx\csname l@#1\endcsname\relax
\typeout{** WARNING: IEEEtran.bst: No hyphenation pattern has been}%
\typeout{** loaded for the language `#1'. Using the pattern for}%
\typeout{** the default language instead.}%
\else
\language=\csname l@#1\endcsname
\fi
#2}}
\providecommand{\BIBdecl}{\relax}
\BIBdecl

\bibitem{asher2013introduction}
K.~Asher, ``An introduction to laplace transform,'' \emph{International Journal
  of Science and Research (IJSR), India}, vol.~2, no.~1, pp. 601--606, 2013.

\bibitem{oldham1974fractional}
K.~Oldham and J.~Spanier, ``The fractional calculus, academic press, new
  york,'' \emph{The fractional calculus. Academic Press, New York.}, 1974.

\bibitem{miller1993introduction}
K.~Miller and B.~Ross, ``An introduction to the fractional integrals and
  derivatives-theory and applications. john willey \& sons,'' \emph{Inc., New
  York}, 1993.

\bibitem{kilbas2006theory}
A.~A. Kilbas, H.~M. Srivastava, and J.~J. Trujillo, \emph{Theory and
  applications of fractional differential equations}.\hskip 1em plus 0.5em
  minus 0.4em\relax elsevier, 2006, vol. 204.

\bibitem{klimek2002lagrangean}
M.~Klimek, ``Lagrangean and hamiltonian fractional sequential mechanics,''
  \emph{Czechoslovak Journal of Physics}, vol.~52, no.~11, pp. 1247--1253,
  2002.

\bibitem{agrawal2002formulation}
O.~P. Agrawal, ``Formulation of euler--lagrange equations for fractional
  variational problems,'' \emph{Journal of Mathematical Analysis and
  Applications}, vol. 272, no.~1, pp. 368--379, 2002.

\bibitem{baleanu2006fractional}
D.~Baleanu and O.~P. Agrawal, ``Fractional hamilton formalism within caputo’s
  derivative,'' \emph{Czechoslovak Journal of Physics}, vol.~56, no. 10-11, pp.
  1087--1092, 2006.

\bibitem{rabei2004potentials}
E.~M. Rabei, T.~S. Alhalholy, and A.~Rousan, ``Potentials of arbitrary forces
  with fractional derivatives,'' \emph{International journal of modern physics
  A}, vol.~19, no. 17n18, pp. 3083--3092, 2004.

\bibitem{rabei2006quantization}
E.~M. Rabei, A.-W. Ajlouni, and H.~B. Ghassib, ``Quantization of brownian
  motion,'' \emph{International Journal of theoretical physics}, vol.~45,
  no.~9, pp. 1613--1623, 2006.

\bibitem{rabei2007hamilton}
E.~M. Rabei, K.~I. Nawafleh, R.~S. Hijjawi, S.~I. Muslih, and D.~Baleanu, ``The
  hamilton formalism with fractional derivatives,'' \emph{Journal of
  Mathematical Analysis and Applications}, vol. 327, no.~2, pp. 891--897, 2007.

\bibitem{rabeihamilton}
E.~M. Rabei and B.~S. Ababneh, ``Hamilton-jacobi fractional sequential
  mechanics.''

\bibitem{hilfer2000applications}
R.~Hilfer \emph{et~al.}, \emph{Applications of fractional calculus in
  physics}.\hskip 1em plus 0.5em minus 0.4em\relax World scientific Singapore,
  2000, vol.~35, no.~12.

\bibitem{podlubny1998fractional}
I.~Podlubny, \emph{Fractional differential equations: an introduction to
  fractional derivatives, fractional differential equations, to methods of
  their solution and some of their applications}.\hskip 1em plus 0.5em minus
  0.4em\relax Elsevier, 1998.

\bibitem{khalil2014new}
R.~Khalil, M.~Al~Horani, A.~Yousef, and M.~Sababheh, ``A new definition of
  fractional derivative,'' \emph{Journal of Computational and Applied
  Mathematics}, vol. 264, pp. 65--70, 2014.

\bibitem{atraoui2021existence}
M.~Atraoui and M.~Bouaouid, ``On the existence of mild solutions for nonlocal
  differential equations of the second order with conformable fractional
  derivative,'' \emph{Advances in Difference Equations}, vol. 2021, no.~1, pp.
  1--11, 2021.

\bibitem{atangana2015new}
A.~Atangana, D.~Baleanu, and A.~Alsaedi, ``New properties of conformable
  derivative,'' \emph{Open Mathematics}, vol.~1, no. open-issue, 2015.

\bibitem{singh2018numerical}
B.~K. Singh and A.~Kumar, ``Numerical study of conformable space and time
  fractional fokker--planck equation via cfdt method,'' in \emph{International
  Conference on Recent Advances in Pure and Applied Mathematics}.\hskip 1em
  plus 0.5em minus 0.4em\relax Springer, 2018, pp. 221--233.

\bibitem{shihab2021associated}
H.~Shihab and T.~Y. Al-khayat, ``Associated conformable fractional legendre
  polynomials,'' in \emph{Journal of Physics: Conference Series}, vol. 1999,
  no.~1.\hskip 1em plus 0.5em minus 0.4em\relax IOP Publishing, 2021, p.
  012091.

\bibitem{khalil2019geometric}
R.~Khalil and M.~A. H. M.~A. Hammad, ``Geometric meaning of conformable
  derivative via fractional cords,'' \emph{J. Math. Computer Sci}, vol.~19, pp.
  241--245, 2019.

\bibitem{hammad2014abel}
M.~A. Hammad and R.~Khalil, ``Abel's formula and wronskian for conformable
  fractional differential equations,'' \emph{International Journal of
  Differential Equations and Applications}, vol.~13, no.~3, 2014.

\bibitem{ahmad2015antiperiodic}
B.~Ahmad, J.~Losada, and J.~J. Nieto, ``On antiperiodic nonlocal three-point
  boundary value problems for nonlinear fractional differential equations,''
  \emph{Discrete Dynamics in Nature and Society}, vol. 2015, 2015.

\bibitem{al2021extension}
M.~Al-Masaeed, E.~M. Rabei, A.~Al-Jamel, and D.~Baleanu, ``Extension of
  perturbation theory to quantum systems with conformable derivative,''
  \emph{Modern Physics Letters A}, p. 2150228, 2021.

\bibitem{chung2021new}
W.~S. Chung, H.~Hassanabadi, and E.~Maghsoodi, ``A new fractional mechanics
  based on fractional addition,'' \emph{Revista Mexicana de F{\'\i}sica},
  vol.~67, no.~1, pp. 68--74, 2021.

\bibitem{AlMasaeedRabeiAlJamelBaleanu+2021+395+401}
\BIBentryALTinterwordspacing
M.~Al-Masaeed, E.~M. Rabei, A.~Al-Jamel, and D.~Baleanu, ``Quantization of
  fractional harmonic oscillator using creation and annihilation operators,''
  \emph{Open Physics}, vol.~19, no.~1, pp. 395--401, 2021. [Online]. Available:
  \url{https://doi.org/10.1515/phys-2021-0035}
\BIBentrySTDinterwordspacing

\bibitem{al2021wkb}
M.~Al-Masaeed, E.~Rabei, A.~Al-Jamel \emph{et~al.}, ``Wkb approximation with
  conformable operator,'' \emph{arXiv preprint arXiv:2111.01547}, 2021.

\bibitem{al2021effect}
A.~Al-Jamel, M.~Al-Masaeed, E.~Rabei, D.~Baleanu \emph{et~al.}, ``The effect of
  deformation of special relativity by conformable derivative,'' \emph{arXiv
  preprint arXiv:2111.02799}, 2021.

\bibitem{mozaffari2018investigation}
F.~Mozaffari, H.~Hassanabadi, H.~Sobhani, and W.~Chung, ``Investigation of the
  dirac equation by using the conformable fractional derivative,''
  \emph{Journal of the Korean Physical Society}, vol.~72, no.~9, pp. 987--990,
  2018.

\bibitem{alextension}
M.~Al-Masaeed, E.~M. Rabei, and A.~Al-Jamel, ``Extension of the variational
  method to conformable quantum mechanics,'' \emph{Mathematical Methods in the
  Applied Sciences}.

\bibitem{hammad2021analytical}
M.~Hammad, A.~S. Yaqut, M.~Abdel-Khalek, and S.~Doma, ``Analytical study of
  conformable fractional bohr hamiltonian with kratzer potential,''
  \emph{Nuclear Physics A}, vol. 1015, p. 122307, 2021.

\bibitem{yavuz2018conformable}
M.~Yavuz and B.~Ya{\c{s}}k{\i}ran, ``Conformable derivative operator in
  modelling neuronal dynamics.'' \emph{Applications \& Applied Mathematics},
  vol.~13, no.~2, 2018.

\bibitem{xin2019modeling}
B.~Xin, W.~Peng, Y.~Kwon, and Y.~Liu, ``Modeling, discretization, and
  hyperchaos detection of conformable derivative approach to a financial system
  with market confidence and ethics risk,'' \emph{Advances in Difference
  Equations}, vol. 2019, no.~1, pp. 1--14, 2019.

\bibitem{kumar2020variety}
D.~Kumar, M.~Kaplan, M.~Haque, M.~Osman, D.~Baleanu \emph{et~al.}, ``A variety
  of novel exact solutions for different models with the conformable derivative
  in shallow water,'' \emph{Frontiers in Physics}, vol.~8, p. 177, 2020.

\bibitem{kumar2018modified}
D.~Kumar, A.~R. Seadawy, and A.~K. Joardar, ``Modified kudryashov method via
  new exact solutions for some conformable fractional differential equations
  arising in mathematical biology,'' \emph{Chinese journal of physics},
  vol.~56, no.~1, pp. 75--85, 2018.

\bibitem{khater2019modified}
M.~Khater, R.~A. Attia, and D.~Lu, ``Modified auxiliary equation method versus
  three nonlinear fractional biological models in present explicit wave
  solutions,'' \emph{Mathematical and Computational Applications}, vol.~24,
  no.~1, p.~1, 2019.

\bibitem{lazo2016variational}
M.~J. Lazo and D.~F. Torres, ``Variational calculus with conformable fractional
  derivatives,'' \emph{IEEE/CAA Journal of Automatica Sinica}, vol.~4, no.~2,
  pp. 340--352, 2016.

\bibitem{abdeljawad2015conformable}
T.~Abdeljawad, ``On conformable fractional calculus,'' \emph{Journal of
  computational and Applied Mathematics}, vol. 279, pp. 57--66, 2015.

\bibitem{boyce2017elementary}
W.~E. Boyce, R.~C. DiPrima, and D.~B. Meade, \emph{Elementary differential
  equations}.\hskip 1em plus 0.5em minus 0.4em\relax John Wiley \& Sons, 2017.

\bibitem{abu2019laguerre}
M.~Abu~Hammad, B.~Albarmawi, A.~Shmasneh, A.~Dababneh \emph{et~al.}, ``Laguerre
  equation and fractional laguerre polynomials,'' \emph{J. Semigroup Theory
  Appl.}, vol. 2019, pp. Article--ID, 2019.

\bibitem{shat2019fractional}
R.~Shat, S.~Alrefai, I.~Alhamayda, A.~Sarhan, and M.~Al-Refai, ``The fractional
  laguerre equation: Series solutions and fractional laguerre functions,''
  \emph{Frontiers in Applied Mathematics and Statistics}, vol.~5, p.~11, 2019.

\bibitem{abdelhakim2019flaw}
A.~A. Abdelhakim, ``The flaw in the conformable calculus: it is conformable
  because it is not fractional,'' \emph{Fractional Calculus and Applied
  Analysis}, vol.~22, no.~2, pp. 242--254, 2019.

\bibitem{al2019new}
Z.~Al-Zhour, F.~Alrawajeh, N.~Al-Mutairi, and R.~Alkhasawneh, ``New results on
  the conformable fractional sumudu transform: theories and applications,''
  \emph{International Journal of Analysis and Applications}, vol.~17, no.~6,
  pp. 1019--1033, 2019.

\bibitem{khader2017conformable}
A.~Khader, ``The conformable laplace transform of the fractional chebyshev and
  legendre polynnomials,'' Ph.D. dissertation, MSc. Thesis, Zarqa University,
  2017.

\end{thebibliography}
\bibliographystyle{IEEEtran}
\end{document}